\documentclass{article}[12pt]
\usepackage[a4paper,  margin=.5in]{geometry}
\usepackage{setspace}
\usepackage{graphicx}
\graphicspath{ {C:/Users/cssc/Desktop/graphs(main)/} }
\usepackage{amsmath}
\usepackage{mathrsfs}
\usepackage{amsfonts}
\usepackage{amssymb}
\usepackage{amsthm}
\usepackage{titlesec}
\usepackage{mathtools}
\usepackage{stackengine}
\usepackage{threeparttable}
\usepackage[dvipsnames,table]{xcolor}
\usepackage[colorlinks,allcolors=red]{hyperref}

\theoremstyle{plain}% Theorem-like structures provided by amsthm.sty
\newtheorem{theorem}{Theorem}[section]

\newtheorem{corollary}[theorem]{Corollary}

\theoremstyle{definition}

\newtheorem{example}[theorem]{Example}

\theoremstyle{remark}

\titleformat{\section}
{\normalfont\Large\bfseries}{\thesection}{1em}{}
\providecommand{\keywords}[1]
{
	\small	
	\textbf{\textit{Keywords---}} #1
}

\usepackage{amsthm}
\newtheorem{manualtheoreminner}{Theorem}
\newenvironment{manualtheorem}[1]{%
	\renewcommand\themanualtheoreminner{#1}%
	\manualtheoreminner
}{\endmanualtheoreminner}
\begin{document}
	\title{\textbf{The extended Bregman divergence and parametric estimation}}

	\author{Sancharee Basak \and Ayanendranath Basu}
	\date{%
		Interdisciplinary Statistical Research Unit\\
		Indian Statistical Institute, Kolkata, INDIA\\%
		\today
	}
	\maketitle

\begin{abstract}
Minimization of suitable statistical distances~(between the data and model densities) has proved to be a very useful technique in the field of robust inference. Apart from the class of $\phi$-divergences of \cite{a} and \cite{b}, the Bregman divergence (\cite{c}) has been extensively used for this purpose. However, since the data density must have a linear presence in the cross product term of the Bregman divergence involving both the data and model densities, several useful divergences cannot be captured by the usual Bregman form. In this respect, we provide an extension of the ordinary Bregman divergence by considering an exponent of the density function as the argument rather than the density function itself. We demonstrate that many useful divergence families, which are not ordinarily Bregman divergences, can be accommodated within this extended description. Using this formulation, one can develop many new families of divergences which may be useful in robust inference. In particular, through an application of this extension, we propose the new class of the GSB divergence family. We explore the applicability of the minimum GSB divergence estimator in discrete parametric models. Simulation studies as well as conforming real data examples are given to demonstrate the performance of the estimator and to substantiate the theory developed.
\end{abstract}

%\begin{keywords}
%\end{keywords}

\section{Introduction}

In the domain of statistical inference, there is, generally, an inherent trade-off between the model  efficiency and the robustness of the procedure.~Often we compromise the efficiency of the procedure to a certain allowable extent for achieving better robustness properties. In the present age of big data, the robustness angle of statistical inference has to be dealt with greater care than ever before. Such methods do exist in the literature which allow full asymptotic efficiency simultaneously with strong robustness properties; see, e.g., \cite{d}, \cite{e} and \cite{f}. In practice, however, one would be hard-pressed to find a procedure which matches the likelihood-based methods in terms of efficiency in small to moderate samples without inheriting any of the robustness limitations of the latter. Many of these trade-off issues are discussed in the canonical texts on robustness such as \cite{g}, \cite{h} and \cite{i}; for a minimum divergence view of this issue, see \cite{f}.

There are several types of divergences which are used in minimum distance inference. Most of them are not mathematical metrics. They may not satisfy the triangle inequality or may not even be symmetric in their arguments. The only properties we demand of these measures are that they are non-negative and are equal to zero if and only if the two arguments are identically equal. Sometimes we will refer to these divergence measures as `statistical distances' or, loosely, as `distances' without any claim to metric properties. 

Most of the density-based divergences in the literature belong to either the class of chi-square type distances~(formally called $\phi$-divergences, $f$-divergences or disparities) or Bregman divergences. See \cite{a}, \cite{b} and \cite{e} for a description of the divergences  of the chi-square type and \cite{c} for Bregman divergences. Although \cite{c} introduced the Bregman divergence in order to be used in convex programming, it has, because of its flexible characteristics, been used in many branches of natural science as well as in areas like information theory and computational geometry. The class of chi-square type distances between two densities $g$ and $f$ includes, for example, the likelihood disparity (LD), the Kullback-Leibler divergence (KLD) and the (twice, squared) Hellinger distance (HD), given by 
\begin{equation}
\label{frst}
{\rm LD}(g, f) = \int g \log\left(\frac{g}{f}\right),~
{\rm KLD}(g, f) =\int f \log\left(\frac{f}{g}\right),~{\rm HD}(g, f) = \frac{1}{2} \int (f^{1/2} - g^{1/2})^2, 
\end{equation}
respectively. Representative members of the class of Bregman divergences include the LD and the squared $L_2$ distance, where 
\begin{equation}
\label{ww}
L_2(g, f) = \int (g - f)^2.
\end{equation}
The LD is the only common member between the class of chi-square type distances and Bregman divergences. 

In the parametric estimation scheme that we consider, the estimator corresponds to the parameter of the model density which is closest to the observed data density in terms of the given divergence, the observed data density being a non-parametric representative of the true unknown density based on the given sample. In case of chi-square type distances, the construction of the data density inevitably requires the use of an appropriate non-parametric smoothing technique, like kernel density estimation, in continuous models (the LD is the only exception). This makes the derivation of the asymptotic properties far more involved, and complicates the computational aspect of this estimation. On the other hand, all the minimum Bregman divergence estimators are M-estimators and, hence, they avoid this density estimation component. 
%The Bregman divergence, which is indexed by some convex function, is used to measure the distance between two points in the mathematical domain. In the statistical domain, it leads to several density-based divergences for different choices of convex functions.

In this paper, our primary aim is to extend the scope of the Bregman divergence by utilizing the powers of densities as arguments, rather than the densities themselves; this leads to the generalized class of the extended Bregman divergences which can then be used to generate new divergences which could provide more refined tools for minimum divergence inference compared to the current state of the art. This is the key idea of this work. Note that the use of the Bregman divergence in statistics is relatively recent; the class of density power divergences by \cite{j}, defined in Section 2, is a prominent example of Bregman divergences having significant applications in statistical inference. Many minimum divergence procedures have natural robustness properties against data contamination and outliers. 
%As our class of divergences become more and more rich and refined we expect that better options for statistical data analysis involving parametric inference will be available. 
The extended Bregman divergence allows us to express several existing divergence families as special cases of it, which is not possible through the ordinary Bregman divergence. Consequently, the extended Bregman idea can be used to generate large super-families of divergences containing, together with the existing divergences, many new and useful divergence families as special cases.

In Section 2, we propose the extended Bregman divergence family. We demonstrate how it allows us to capture well known divergences that are not within the ordinary Bregman class and give a potential route for constructing new divergences. In Section 3, by  considering a specific form of the convex function along with a particular choice of exponent of densities, we construct a large super-family of divergences within the extended Bregman family. Several known divergence families are obtained as special cases of this super-family. Section 4 introduces the corresponding minimum distance estimator, while Section 5 studies its asymptotic properties. Section 6 explores the robustness properties of the estimator based on its influence function. A large scale simulation study is taken up in Section 7, and a tuning parameter selection strategy is discussed in Section 8. The final section has some concluding remarks. 

Before concluding this section, we summarize what we believe to be the main achievements in this paper.
\begin{enumerate}
	\item We provide a simple extension of the Bregman divergence by considering powers of densities~ (instead of the densitites themselves) as arguments. Many divergence families~(which are ordinarily not members of the class of Bregman divergences) can now be looked upon as members of this extended family, and the properties of the corresponding minimum distance estimators may be obtained from the general properties common to all the members of the extended family. 
	\item Ordinarily, all minimum Bregman divergence estimators are also M-estimators. But, through this extension, several minimum distance estimators which are not M-estimators also become a part of this extended family. 
	\item Consideration of exponentiated arguments with a specific choice of the convex function introduces a generalized super-family which we refer to as the GSB divergence family. The power divergence family of \cite{k}, the density power divergence (DPD) family of \cite{j}, the Bregman exponential divergence (BED) of \cite{l}  and the $S$-divergence family of \cite{m} can all be brought under the umbrella of this super-family.
	\item This GSB divergence family consists of three tuning parameters $\alpha$, $\beta$ and $\lambda$. By simultaneously varying all these three parameters, we can generate new divergences (and hence new minimum divergence estimators) which are outside the union of the BED and S-divergence family, but can potentially provide improved performance over both of these classes.
\end{enumerate}

\section{The extended Bregman divergence and special cases}

Being motivated by the problem of convex programming, \cite{c} introduced the Bregman divergence, a measure of dissimilarity between any two vectors in the Euclidean space. In $\mathbb{R}^p$, it has the form
\begin{eqnarray}
\label{obd}
& D_\psi \left(\boldsymbol{x},\boldsymbol{y}\right)= \Bigg\{ \psi\left(\boldsymbol{x}\right)-\psi\left(\boldsymbol{y}\right)-\langle\nabla \psi\left(\boldsymbol{y}\right),\boldsymbol{x}-\boldsymbol{y}\rangle \Bigg\},
\end{eqnarray}
for any strictly convex function $\psi:\mathcal{S}\rightarrow\mathbb{R}$ and for any two $p$-dimensional vectors $\boldsymbol{x},\boldsymbol{y} \in \mathcal{S}$, where $\mathcal{S}$ is a convex subset of $\mathbb{R}^p$. Here, $\nabla \psi\left(\boldsymbol{y}\right)$ denotes the gradient of $\psi$ with respect to its argument at $\boldsymbol{y}=\left(y_1,y_2,\ldots,y_p\right)^{T}$. It is evident that only the convexity criterion of the function $\psi\left(\cdot\right)$ is necessary for the non-negativity property of the divergence  $D_\psi\left(\boldsymbol{x},\boldsymbol{y}\right)$ to hold. One could, therefore, consider other quantities as the arguments rather than the points themselves in this measure. Hence, as long as $\psi$ remains convex, any set of arguments whose equivalence translates to the equivalence of $\boldsymbol{x}$ and $\boldsymbol{y}$ can be used in the distance expression. This observation may be used to extend the Bregman divergence to have the form
\begin{eqnarray}
\label{bd}
& D_\psi\left(\boldsymbol{x},\boldsymbol{y}\right)= \Bigg\{\psi\left(\boldsymbol{x}^k\right)-\psi\left(\boldsymbol{y}^k\right)-\langle\nabla \psi\left(\boldsymbol{y}^k\right),\boldsymbol{x}^k-\boldsymbol{y}^k\rangle\Bigg\}.
\end{eqnarray}
Here, $\nabla \psi\left(\boldsymbol{y}^k\right)$ denotes the gradient of $\psi$ with respect to its argument, evaluated at $\boldsymbol{y}^k=\left(y_1^k,y_2^k,\ldots,y_d^k\right)^{T}$ and $\psi$ is a strictly convex function, mapping $\mathcal{S}$ to $\mathbb{R}$, $\mathcal{S}$ being a convex subset of $\mathbb{R^+}^p$. Since our main purpose is to utilize this extension in the field of statistics where the arguments, being probability density functions, are inherently non-negative, restricting the domain of $\psi$ to $\mathbb{R^+}^p$ does not cause any difficulty. It is also not difficult to see that many of the properties of the Bregman divergence in Equation (\ref{obd}), as described by \cite{n}, are retained by the extended version in Equation  (\ref{bd}). However, we will not make use of these properties in this paper, so we do not discuss them here any further. 

The Bregman divergence has significant applications in the domain of statistical inference for both discrete and continuous models. Given two densities $g$ and $f$, the Bregman divergence between these densities (associated with the convex function $\psi$) is given by 
\begin{eqnarray}
\label{bdd}
& D_\psi \left(g,f\right)=\displaystyle\int\Bigg\{\psi\left(g\left(x\right)\right)-\psi\left(f\left(x\right)\right)-\left(g\left(x\right)-f\left(x\right)\right)\nabla \psi\left(f\left(x\right)\right)\Bigg\}dx.
\end{eqnarray}
By the strict convexity of the function $\psi$, the integrand in the above Equation (\ref{bdd}) 
is non-negative and, therefore, so is the integral. It is also clear that the divergence equals zero if and only if the arguments $g$ and $f$ are identically equal. Well-known examples include the LD and the (squared) $L_{2}$ distance which correspond to $\psi\left(x\right)=x \log x$ and $\psi\left(x\right)=x^2$ respectively. In a real scenario, one uses $g$ as the true data generating density and $f=f_\theta$ as the parametric model density.

In \cite{j}, the important class of density power divergences (DPDs) has been proposed, which is a subfamily of the class of Bregman divergences. This family is generated by the function $\psi\left(x\right)=\frac{x^{\alpha+1}-x}{\alpha}$, indexed by a non-negative tuning parameter $\alpha$. As a function of $\alpha$, the density power divergence may be expressed as 
\begin{eqnarray}
& DPD_{\alpha} \left(g,f\right)=\displaystyle\int\left\{f^{\alpha+1}\left(x\right)-\left(1+\frac{1}{\alpha}\right)g\left(x\right)f^\alpha\left(x\right)+\frac{1}{\alpha}g^{\alpha+1}\left(x\right)\right\}dx.
\end{eqnarray}
It is a simple matter to check that for $\alpha=1$, the above reduces to the (squared) $L_{2}$ distance between $g$ and $f$, whereas when  $\alpha \rightarrow 0$, one recovers the likelihood disparity defined in Equation (\ref{frst}).
% into the area of robust estimation based on minimum dist ${\rm LD}(g,f)$ as given in Equation (\ref{ww}).

The DPD class should not be confused with the power divergence~(PD) class of divergences~(see \cite{k}) which has the form
\begin{eqnarray}
\label{pd}
PD_{\lambda}\left(g,f\right)=\displaystyle \frac{1}{\lambda\left(\lambda+1\right)}\int\left\{g\left(x\right)\left(\frac{g\left(x\right)}{f\left(x\right)}\right)^\lambda-1\right\}dx, \lambda \in \mathbb{R}.
\end{eqnarray}
The PD class is a subfamily of chi-square type distances. The latter class of divergences has the form 
\begin{equation}
\rho\left(g,f\right)=\displaystyle \int C\left(\delta\left(x\right)\right) f\left(x\right) dx,
\end{equation}
where $C$ is a strictly convex function and $\delta\left(x\right)=\frac{g\left(x\right)}{f\left(x\right)}-1$. The power divergence corresponds to the specific convex function 
\begin{equation}
C\left(\delta\right)=\frac{\left(\delta+1\right)^{\lambda+1}-\left(\delta+1\right)}{\lambda\left(\lambda+1\right)}-\frac{\delta}{\lambda+1}.
\end{equation}
Important special cases of the PD class include the LD~(obtained in the limit as $\lambda \rightarrow 0$) and the twice, squared HD~(obtained for $\lambda=-\frac{1}{2}$). The LD is the only common member between the PD and the DPD classes. 

The Bregman exponential divergence (BED) class~(\cite{l}), on the other hand, has the form 
\begin{eqnarray}
\label{bb}
& BED_\beta\left(g,f\right)=\displaystyle\frac{2}{\beta}\displaystyle\int\left\{e^{\beta f\left(x\right)}\left(f\left(x\right)-\frac{1}{\beta}\right)-e^{\beta f\left(x\right)}g\left(x\right)+\frac{e^{\beta g\left(x\right)}}{\beta}\right\}dx.
\end{eqnarray}
The defining convex function is $\psi\left(x\right)=\frac{2\left(e^{\beta x}-\beta x-1\right)}{\beta^2}$ which is indexed by the real parameter $\beta$. This family generates the (squared) $L_{2}$ distance in the limit $\beta \rightarrow 0$. 

A list of some Bregman divergences useful in the context of statistical inference is presented in Table \ref{tab1}.

\begin{table}[h]
	\begin{center}
		\caption{Different divergences as special cases of the Bregman divergence}
		\begin{threeparttable}
			\label{tab1}
			\renewcommand{\arraystretch}{2.5}
			\begin{tabular}{c c}
				\hline
				%\multicolumn{4}{c}{Different divergences as special cases of Bregman Divergence}\\
				\hline
				Choice of convex function & Divergences\\
				\hline
				$B\left(x\right)=\displaystyle x^2$ & (squared) $L_{2}$ Distance\\
				$B\left(x\right)=\displaystyle x \log \left(x\right)$ & Likelihood Disparity\\
				$B\left(x\right)=\displaystyle \frac{x^{1+\alpha}-x}{\alpha}$ & Density Power Divergence~(DPD)\\
				$B\left(x\right)=\displaystyle -\frac{\log \left(x\right)}{2\pi}$ & Itakura-Saito Distance\\
				$B\left(x\right)=\frac{2\left(e^{\beta x}-\beta x-1\right)}{\beta^2}$ & Bregman Exponential Divergence\\
				\hline
			\end{tabular}
		\end{threeparttable}
	\end{center}
\end{table}

Consider the standard set up of parametric estimation where $G$ is the true data generating distribution which is modeled by the parametric family 
${\cal F} = \{F_\theta: \theta \in \Theta \subset {\mathbb R}^p\}$. Let $g$ and $f_\theta$ be the corresponding density functions. Further we assume that both $G$ and $F_{\theta}$ belong to $\mathcal{G}$, the class of all cumulative distribution functions having densities with respect to some appropriate dominating measure. Our aim is to estimate the unknown parameter $\theta$ by choosing the model density closest to the true density in the Bregman sense. The definition of ordinary Bregman divergences as given in Equation (\ref{bdd}), useful as it is, does not include many well-known and popular divergences which are extensively used in the literature for different purposes including parameter estimation. The PD family is a prominent example. An inspection of the Bregman form in Equation (\ref{bdd}) indicates that the term which involves both densities $g$ and $f$ is of the form 
\begin{equation}
\label{se}
\int g\left(x\right)\nabla \psi\left(f\left(x\right)\right)dx.
\end{equation}
Here, the density $g$ is present only as a linear term having exponent one. Given a random sample $X_1, X_2, \ldots, X_n$ from the true distribution $G$, the term in Equation (\ref{se})  can be empirically estimated by $\frac{1}{n} \sum \nabla \psi(f_\theta(X_i))$ (with $f = f_\theta$ under the parametric model) so that one can construct an empirical version of the divergence without any non-parametric density estimation. On the other hand, this restricts the class of divergences that are expressible in the Bregman form. Using an extension in the spirit of Equation (\ref{bd}) may allow the construction of richer classes of divergences. With this aim, we define the \textit{extended Bregman divergence} between two densities $g$ and $f$ as 
\begin{eqnarray}
\fontsize{7.5pt}{8}\selectfont
\label{ss1}
& D^{(k)}_\psi \left(g,f\right)=\displaystyle \int\left\{\psi\left(g^k\left(x\right)\right)-\psi\left(f^k\left(x\right)\right)-\left(g^k\left(x\right)-f^k\left(x\right)\right)\nabla \psi\left(f^k\left(x\right)\right)\right\}dx.\nonumber\\
\end{eqnarray}
Apart from the requirement of strict convexity of the  function $\psi$, this formulation also depends on a positive index $k$ with which the density is exponentiated. For the rest of the paper, the notation $D^{(k)}_\psi(\cdot, \cdot)$ will refer to this general form in Equation (\ref{ss1}), of which the divergence in Equation (\ref{bdd}) is a special case for $k =1$. Evidently, $D^{(k)}_\psi\left(g,f\right)\geq0$ for any choices of densities $f$ and $g$ with respect to the same measure. Moreover, the fact that $D^{(k)}_\psi\left(g,f\right)=0$ if and only if $g=f$, holds true in this case due to non-negativity property of a density as well as the consideration of strict convexity of the function $\psi\left(\cdot\right)$.

In the following, we present some special cases of extended Bregman divergence.

\begin{enumerate}
	\item \textbf{$S$-Hellinger family~(\cite{m})}
	If we take $\psi\left(x\right)=\frac{2e^{\beta x}}{\beta^2}$ with $k=\frac{1+\alpha}{2}$, $\alpha \in \left(0,1\right)$ in Equation (\ref{ss1}), it will generate an extension of the BED family having the form
	\begin{equation}
	\label{ee}
	BED^{(k)}_{\beta}\left(g,f\right)=\frac{2}{\beta}\int\left\{e^{\beta f^k\left(x\right)}\left(f^k\left(x\right)-\frac{1}{\beta}\right)-e^{\beta f^k\left(x\right)}g^k\left(x\right)+\frac{e^{\beta g^k\left(x\right)}}{\beta}\right\}dx.
	\end{equation}
	It can be easily shown that, as $\beta \rightarrow 0$ and $k=\frac{1+\alpha}{2}$, $\alpha \in \left(0,1\right)$, the application of L'Hospital's rule leads to the S-Hellinger Distance (SHD) family with the form
	\begin{equation}
	SHD_\alpha\left(g,f\right)=\frac{2}{1+\alpha}\int\left(g^{\frac{1+\alpha}{2}}\left(x\right)-f^{\frac{1+\alpha}{2}}\left(x\right)\right)^2 dx.
	\end{equation}
	This was introduced by \cite{m} as a special case of the $S$-divergence family. This family cannot be expressed through the normal expression of the Bregman divergence, but through this extension, we can express this member of the $S$-divergence family as a~(limiting) member of the extended BED class.
	\item \textbf{PD family~(\cite{k})}
	If we take $\psi\left(x\right)=\frac{x^{1+\frac{B}{A}}}{B}$, $A=1+\lambda$, $B=-\lambda$ and $\lambda \in \mathbb{R}$, with $k=A$ in Equation (\ref{ss1}), we get the PD family introduced in Equation (\ref{pd}).
	\item \textbf{$S$-divergence family~(\cite{m})}
	If we take $\psi\left(x\right)=\frac{x^{1+\frac{B}{A}}}{B}$, $A=1+\lambda\left(1-\alpha\right)$, $B=\alpha-\lambda\left(1-\alpha\right)$, $A+B=1+\alpha$, $\alpha\geq0$, $\lambda\in \mathbb{R}$ and $k=A$ in Equation (\ref{ss1}), we get the $S$-divergence having the following form
	\begin{equation}
	\fontsize{9pt}{8}\selectfont
	SD_{(\alpha, \lambda)}\left(g,f\right) = \int\left\{\frac{1}{B}\left(g^{A+B}\left(x\right)-f^{A+B}\left(x\right)\right)-\left(g^A\left(x\right)-f^A\left(x\right)\right)\frac{A+B}{AB}f^B\left(x\right)\right\}dx.
	\end{equation}
	This is one of the most useful divergence families in the domain of robust inference due to its capacity to generate much more robust estimator(s) than the DPD and PD families can generate. Through this extension of Bregman divergence, it is now possible to express this divergence as a special case of the extended family.
	%\item \textbf{BED family~(Mukherjee et al.~(2019))}
	%We have already shown the derivation of the extended BED family through our proposal in Equation (\ref{ee}). More specifically if we choose $A=1$ here, then we get the normal BED family given in Equation (\ref{bb}).
	%Some simple algebra shows that it can be written in the form $\int C\left(\delta\left(x\right)\right)f\left(x\right)dx$ where $C\left(\delta\right)$ is exactly equal to the disparity generating function of the PD family.
\end{enumerate}

Note that, for $k \neq 1$, $f^k$ and $g^k$ will generally no longer represent probability densities, and by extending the divergence idea to general positive measures (beyond probability measures), \cite{o} has suggested certain constructions where the power divergence has been exhibited in the Bregman divergence form for general measures. Also, see the discussion in \cite{p}. We differ from the interpretation in these papers in the sense that we still view the family of divergences presented here in Equation (\ref{ss1}) as divergences between valid probability densities. Given any two probability densities, the expression in Equation (\ref{ss1}) is non-negative, and equals zero if and only if the densities $g$ and $f$ are identical, irrespective of the value of $k$. 

For the ordinary Bregman divergence, the term in Equation (\ref{se}), with $f = f_{\theta}$,  may be approximated by 
\begin{equation}
\frac{1}{n}\sum_{i=1}^{n} \nabla \psi\left(f_{\theta}\left(X_{i}\right)\right),
\end{equation} through the replacement of $dG\left(x\right)$ by $dG_n\left(x\right)$, where $G_n$ is the empirical distribution function obtained from the random sample $X_1, X_2, \ldots, X_n$. It is evident that the minimizer of the empirical version of the Bregman divergence is an M-estimator. But there are several useful divergences (or divergence families) where the empirical representation of the term involving $f$ and $g$ is not possible using the above trick, and such divergences generate estimators beyond the M-estimator class. See \cite{q} for more discussion on this issue. In the above we have given several examples where the extended Bregman class contains such divergences which are not covered  by the ordinary Bregman form. Thus the structure of the extended class allows us to extend the scope much beyond that of the ordinary Bregman divergence.

\section{Introducing a new divergence family}
Our aim here is to exploit the extended Bregman idea and generate rich new super families of divergences by choosing a suitable convex generating function and a suitable exponent. In particular, we use the convex function $\psi(x)=e^{\beta x}+\frac{x^{1+\frac{B}{A}}}{B}$, $A=1+\lambda\left(1-\alpha\right)$, $B=\alpha-\lambda\left(1-\alpha\right)$, $A+B=1+\alpha$, $\alpha \geq -1$, $\beta, \lambda \in \mathbb{R}$, which, together with the exponent $k=A$, generates the divergence
\begin{equation}
\fontsize{9pt}{8}\selectfont
\label{ebg}
D^*\left(g,f\right) = \int\left\{e^{\beta f^A}\left(\beta f^A-\beta g^A-1\right)+e^{\beta g^A}+\frac{1}{B}\left(g^{A+B}-f^{A+B}\right)-\left(g^A-f^A\right)\frac{A+B}{AB}f^B\right\}dx,
\end{equation}
which we refer to as the GSB divergence (being the abbreviated form of generalized S-Bregman divergence). The divergence measure $D^*$ is also a function of $\alpha, \lambda$ and $\beta$, which we suppress for brevity. 

If we put $A+B=0$ in the above expression with $A\neq0$ and $B\neq0$, we will get the extended BED family with parameter $\beta$ and exponent $k=A$. Moreover, if $A=1$, i.e., $\lambda=0$ then it will lead to the ordinary BED family with parameter $\beta$. On the contrary, if we put $\beta=0$, it will lead to the $S$-divergence family with parameters $\alpha$ and $\lambda$ (in terms of $A$ and $B$). More specifically, when $\alpha=0$ and $\beta=0$, it leads to the power divergence (PD) family. On the other hand, $\beta=0$ and $\lambda=0$ lead us to the density power divergence~(DPD) family. Thus, it acts as a connector between the BED and the $S$-divergence family.

\subsection{Special cases}
We will get several well-known divergences or divergence families from the general form of GSB for particular choices of the three tuning parameters $\alpha$, $\lambda$ and $\beta$. Some such choices are given in Table \ref{tab2}.

\begin{table}[h]
	\begin{center}
		\caption{Different divergences as special cases of GSB divergence}
		\begin{threeparttable}
			\label{tab2}
			\renewcommand{\arraystretch}{1.2}
			\begin{tabular}{c c c c}
				\hline
				$\alpha$ & $\lambda$ & $\beta$ & Divergences\\
				\hline
				$\alpha=-1$ & $\lambda=0$ & $\beta \in \mathbb{R}$ & Bregman Exponential Divergence\tnote{1}\\
				$\alpha=0$ & $\lambda \in \mathbb{R}$ & $\beta=0$ & Power Divergence\\
				$\alpha \in \mathbb{R}$ & $\lambda=0$ & $\beta=0$ & Density Power Divergence\\
				$\alpha \in \mathbb {R}$ & $\lambda \in \mathbb{R}$ & $\beta=0 $ & $S$-Divergence\\
				$\alpha=0$ & $\lambda=-1$ & $\beta=0$ & Kullback-Liebler Divergence\\
				$\alpha=0$ & $\lambda=0$ & $\beta=0$ & Likelihood Disparity\\
				$\alpha=0$ & $\lambda=-.5$ & $\beta=0$ & Hellinger Distance\\
				$\alpha \in \mathbb{R}$ & $\lambda=-.5$ & $\beta=0$ & S-Hellinger Distance\\
				$\alpha=0$ & $\lambda=1$ & $\beta=0$ & Pearson's Chi-square Divergence\\
				$\alpha=0$ & $\lambda=-2$ & $\beta=0$ & Neyman's Chi-square Divergence\\
				$\alpha=1$ & $\lambda \in \mathbb{R}$ & $\beta=0$ & (squared) $L_2$ Distance\\
				\hline
			\end{tabular}
			\begin{tablenotes}
				\footnotesize
				\item[1] This is a constant time B-exponential divergence. It basically generates all the members of BED family corresponding to the same $\beta$ except (squared) $L_{2}$ distance, which occurs when $\beta \rightarrow 0$. However, as seen above, the (squared) $L_2$ distance remains a member of the GSB class for other choices of the tuning parameters.
			\end{tablenotes}
		\end{threeparttable}
	\end{center}
\end{table}

\section{The minimum GSB divergence estimator}
Under the parametric set-up described in Section 2, we would like to identify the best fitting parameter $\theta^g$ by choosing the element of the model family of distributions which provides the closest match to the true density $g$ in terms of the given divergence. The minimum GSB divergence functional $T_{\alpha,\lambda,\beta}:\mathcal{G}\rightarrow\Theta$ is defined by the relation
$$D^*\left(g,f_{T_{\alpha,\lambda,\beta}}\right)=\min\{D^*\left(g,f_\theta\right):\theta \in \Theta\},$$ provided the minimum exists. If the parametric model family is identifiable, it follows from the definition of the divergence that $D^*\left(g,f_\theta\right)=0$, if and only if $g=f_\theta$. Thus, $T_{\alpha,\lambda,\beta}\left(F_\theta\right)=\theta$, uniquely. Hence, the functional $T_{\alpha,\lambda,\beta}$ is Fisher consistent. Given the density $g$, a straightforward differentiation of the GSB divergence of Equation (\ref{ebg}) leads to the estimating equation 
\begin{eqnarray}
\label{a}
\displaystyle \int \left\{A^2 \beta^2 e^{\beta f_{\theta}^A\left(x\right)}f_{\theta}^A\left(x\right) +\left(A+B\right)f_{\theta}^{B}\left(x\right)\right\}\left(f_{\theta}^A\left(x\right)-g^A\left(x\right)\right)u_{\theta}\left(x\right)dx=0.
\end{eqnarray}
In practice, the true density $g$ is unknown, so one has to use a suitable non-parametric density estimator $\hat{g}$ for $g$, depending on the situation. Under a discrete parametric set-up, the natural choice for $\hat{g}$ is the vector of relative frequencies as obtained from the sample data. Thus, assuming a discrete parametric model, and assuming, without loss of generality, that the support of the random variable is $\{0, 1, 2, \ldots,\}$, the estimating equation becomes 
\begin{eqnarray}
\label{ba}
%&A^2 \beta^2\int\{e^{\beta f_{\theta}^A\left(x\right)}f_{\theta}^A\left(x\right) u_\theta\left(x\right)\left(f_{\theta}^A\left(x\right)-\hat{g}^A\left(x\right)\right)\}u_{\theta}\left(x\right)dx+\left(A+B\right)\int\{f_{\theta}^{A+B}\left(x\right)-f_{\theta}^B\left(x\right)\hat{g}^A\left(x\right)\}u_{\theta}\left(x\right)dx=0\nonumber\\
&&\displaystyle \sum_{x=0}^{\infty} A^2 \beta^2e^{\beta f_{\theta}^A\left(x\right)}f_{\theta}^{2A}\left(x\right) u_\theta\left(x\right)+\sum_{x=0}^{\infty} \left(A+B\right)f_{\theta}^{A+B}\left(x\right) u_\theta\left(x\right)\nonumber\\&=&\displaystyle\sum_{x=0}^{\infty}A^2 \beta^2 e^{\beta f_{\theta}^A\left(x\right)}f_{\theta}^A\left(x\right) \hat{g}^A\left(x\right)u_{\theta}\left(x\right)+\displaystyle \sum_{x=0}^{\infty}\left(A+B\right) f_{\theta}^B\left(x\right)\hat{g}^A\left(x\right)u_{\theta}\left(x\right).
\end{eqnarray}
For continuous models, on the other hand, some suitable non-parametric  smoothing technique such as kernel density estimation is inevitable unless the exponent $A$ equals 1.  In the latter case, $g$ appears as a linear term in the estimating equation (\ref{a}). In that case, we can use $G_n$, the empirical distribution function, as an estimator of $G$. Hence, for $A=1$, Equation (\ref{ba}) can be reduced as
\begin{eqnarray}
&&\displaystyle \sum_{x=0}^{\infty}\beta^2e^{\beta f_{\theta}\left(x\right)}f_{\theta}^2\left(x\right) u_\theta\left(x\right)+\sum_{x=0}^{\infty}\left(1+B\right)f_{\theta}^{1+B}\left(x\right) u_\theta\left(x\right) \nonumber\\
&=&\frac{1}{n}\displaystyle \sum_{i=1}^{n}\beta^2 e^{\beta f_{\theta}\left(X_i\right)} f_{\theta}\left(X_i\right) u_\theta\left(X_i\right)+\frac{1}{n}\sum_{i=1}^{n}\left(1+B\right)f_{\theta}^B\left(X_i\right)u_\theta\left(X_i\right).
\end{eqnarray}
Since the left hand side of the above equation is non-random and the right hand side is a sum of independent and identically distributed terms, it is of the form $\sum_{i=1}^{n} \psi\left(X_i,\theta\right)=0$ and the corresponding estimator belongs to the M-estimator class.

In accordance with the information on the first three rows of Table \ref{tab2}, we will refer to the parameters $\alpha$, $\lambda$ and $\beta$ as the DPD parameter, the PD parameter and the BED parameter, respectively. 

\section{Asymptotic properties of the GSB divergence}

In this section, we concentrate on the asymptotic properties of our proposed minimum divergence estimator. As mentioned, we will focus on the discrete set-up throughout the rest of the paper. Let $X_1,X_2,\ldots,X_n$ be independent and identically distributed observations from an unknown distribution $G$ with support $\chi=\{0,1,2,3,\ldots,\}$. On the other hand, we consider a parametric family of distributions $\mathcal{F}=\{F_\theta:\theta \in \Theta \subseteq \mathbb{R}^p \}$, also supported on $\chi$, to  model the true data generating distribution $G$. In this set-up, we assume both $G$ and $\mathcal{F}$ to have densities $g$ and $f_{\theta}$ with respect to the counting measure. Let $\theta^g = T_{\alpha, \beta, \lambda}(G)$ be the best fitting parameter. Since $G$ is unknown, we are going to use the vector of relative frequencies, obtained from the data, as an estimate of $g$ throughout the rest of this paper. Let $r_n(x)$ be the relative frequency of the value $x$ in the sample. The minimum GSB divergence estimator is obtained as a root of the estimating equation
\begin{eqnarray}
&&\displaystyle \sum_{x=0}^{\infty} \left\{A^2 \beta^2 e^{\beta f_{\theta}^A\left(x\right)}f_{\theta}^A\left(x\right) +\left(A+B\right)f_{\theta}^{B}\left(x\right)\right\}\left(f_{\theta}^A\left(x\right)-\hat{g}^A\left(x\right)\right)u_{\theta}\left(x\right)=0\nonumber\\
 &\Rightarrow& \displaystyle \sum_{x=0}^{\infty} \left\{A^2 \beta^2 e^{\beta f_{\theta}^A\left(x\right)}f_{\theta}^{2A}\left(x\right) +\left(A+B\right)f_{\theta}^{A+B}\left(x\right)\right\}\frac{\left(1-\frac{\hat{g}^A\left(x\right)}{f_{\theta}^A\left(x\right)}\right)}{A}u_{\theta}\left(x\right)=0\nonumber\\
 &\Rightarrow& \displaystyle \sum_{x=0}^{\infty} K\left(\delta\left(x\right)\right)\left(A^2 \beta^2 e^{\beta f_{\theta}^A\left(x\right)}f_{\theta}^{2A}\left(x\right)+\left(A+B\right) f_{\theta}^{A+B}\left(x\right)\right)u_\theta \left(x\right)=0,
\end{eqnarray}
where, $\delta\left(x\right)=\delta_n\left(x\right)=\frac{\hat{g}\left(x\right)}{f_\theta\left(x\right)}-1=\frac{r_n\left(x\right)}{f_\theta\left(x\right)}-1$, $K\left(\delta\right)=\frac{\left(\delta+1\right)^A-1}{A}$ and $u_\theta\left(x\right)$ is the likelihood score function at $x$. We denote the minimum GSB divergence estimator obtained as a solution of the above equation as $\hat{\theta}$. Let 
\begin{eqnarray}
J_g&=&J_{\alpha,\beta,A}\left(g\right)\nonumber\\&=&\displaystyle E_g\left(u_{\theta^g}\left(X\right)u^{T}_{\theta^g}\left(X\right)K'\left(\delta_g^g\left(X\right)\right)\left(\left(A+B\right)f_{\theta^g}^\alpha\left(X\right)+A^2\beta^2e^{\beta f_{\theta^g}^A\left(X\right)}f_{\theta^g}^{2A-1}\left(X\right)\right)\right)\nonumber\\&+&\sum_{x=0}^{\infty} K\left(\delta_g^g\left(x\right)\right)\left(\left(A+B\right)f_{\theta^g}^{1+\alpha}\left(x\right)+A^2\beta^2e^{\beta f_{\theta^g}^A\left(x\right)}f_{\theta^g}^{2A}\left(x\right)\right)i_{\theta^g}\left(x\right)\nonumber\\&-&\sum_{x=0}^{\infty} K\left(\delta_g^g\left(x\right)\right)\left(\left(A+B\right)^2 f_{\theta^g}^{1+\alpha}\left(x\right)+A^3\beta^2e^{\beta f_{\theta^g}^A\left(x\right)}f_{\theta^g}^{2A}\left(x\right)\left(2+\beta f_{\theta^g}^A\left(x\right)\right)\right)u_{\theta^g}\left(x\right)u^{T}_{\theta^g}\left(x\right)\nonumber\\
V_g&=&\displaystyle Var_{g}\left(u_{\theta^g}\left(X\right)K'\left(\delta_g^g\left(X\right)\right)\left(\left(A+B\right)f_{\theta^g}^\alpha\left(X\right)+A^2\beta^2e^{\beta f_{\theta^g}^A\left(X\right)}f_{\theta^g}^{2A-1}\left(X\right)\right)\right),
\end{eqnarray}
where, $X$ is a random variable having density $g$, $Var_g$ represents variance under the density $g$, $\delta_g\left(x\right)=\frac{g\left(x\right)}{f_\theta\left(x\right)}-1$, $K'\left(\cdot\right)$ is the derivative of $K\left(\cdot\right)$ with respect to its argument, $\delta_g^g\left(x\right)=\frac{g\left(x\right)}{f_\theta^g\left(x\right)}-1$ and $i_{\theta}\left(x\right)=-u'_{\theta}\left(x\right)$, the negative of the derivative of the score function with respect to the parameter.

\begin{manualtheorem}{1}
	Under the above-mentioned set-up and certain regularity assumptions given in the Online Supplement, there exists a consistent sequence of roots $\hat{\theta}_n$ of the estimating equation (\ref{a}). Moreover, the asymptotic distribution of $\sqrt{n}\left(\hat{\theta}_{n}-\theta^g\right)$ is p-dimensional normal with mean $0$ and $J_g^{-1}V_gJ_g^{-1}$.
\end{manualtheorem}

\begin{corollary}
	When $g=f_\theta$ for some $\theta \in \Theta$, then $\sqrt{n}\left(\theta_n-\theta\right)\sim N\left(0, J^{-1}VJ^{-1}\right)$ asymptotically, where,
	\begin{eqnarray}
	J&=&E_{f_\theta}\left \{u_\theta \left(X\right)u^{T}_\theta \left(X\right)\left(\left(A+B\right)f^\alpha\left(x\right)+A^2 \beta^2 e^{\beta f_{\theta}^A\left(X\right)}f_{\theta}^{2A-1}\left(X\right)\right)\right\}\nonumber\\&=&\sum_{x=0}^{\infty}\left\{u_\theta \left(x\right)u^{T}_\theta \left(x\right)\left(\left(A+B\right)f^\alpha\left(x\right)+A^2 \beta^2 e^{\beta f_{\theta}^A\left(x\right)}f_{\theta}^{2A-1}\left(x\right)\right)\right\}f_\theta\left(x\right),\nonumber\\
	\end{eqnarray}
	
	\begin{eqnarray}
	V& = &V_{f_\theta}\left\{u_\theta \left(X\right)\left(\left(A+B\right)f_\theta^\alpha\left(X\right)+A^2 \beta^2 e^{\beta f_{\theta}^A\left(X\right)}f_{\theta}^{2A-1}\left(X\right)\right)\right\}\nonumber\\
	& = & \left(A+B\right)^2 \sum_{x=0}^{\infty} f^{1+2\alpha}\left(x\right)u_\theta \left(x\right)u^{T}_\theta \left(x\right)\nonumber\\
	&+&A^4 \beta^4 \sum_{x=0}^{\infty} e^{2\beta f_{\theta}^A\left(x\right)}f_{\theta}^{4A-1}\left(x\right)u_\theta \left(x\right)u^{T}_\theta \left(x\right)\nonumber\\
	&+&2\left(A+B\right) A^2 \beta^2 \sum_{x=0}^{\infty} e^{\beta f_{\theta}^A\left(x\right)}f_{\theta}^{2A+\alpha}\left(x\right)u_\theta \left(x\right)u^{T}_\theta \left(x\right)-\zeta\zeta',
	\end{eqnarray}
	where, $\zeta=\displaystyle\sum_{x=0}^{\infty} u_\theta \left(x\right)\left(\left(A+B\right)f^{A+B}\left(x\right)+A^2 \beta^2 e^{\beta f_{\theta}^A\left(x\right)}f_{\theta}^{2A}\left(x\right)\right).$
\end{corollary}

\section{Influence analysis of the minimum GSB estimator}
Here we study the stability of our proposed class of estimators on the basis of the influence function (IF), which measures the effect of adding an infinitesimal mass to the distribution and is one of the most important heuristic tools of robustness. A simple differentiation of a contaminated version of the estimating equation (\ref{a}) leads to the expression
\begin{equation}
IF\left(y,G,T_{\alpha,\lambda,\beta}\right)=J_{G}^{-1}N_{G}\left(y\right),~~\rm{where},
\end{equation}
\begin{eqnarray}
N_{G}\left(y\right)&=&\displaystyle \left(A^2 \beta^2 e^{\beta f^A_{\theta^g}\left(y\right)}f^{A}_{\theta^g}\left(y\right)+(A+B)f^B_{\theta^g}\left(y\right)\right)g^{A-1}\left(y\right)u_{\theta^g}\left(y\right)\nonumber\\&-&\sum_{x=0}^{\infty}\left(A^2 \beta^2 e^{\beta f^A_{\theta^g}\left(x\right)}f^{A}_{\theta^g}\left(x\right)+(A+B)f^B_{\theta^g}\left(x\right)\right)g^{A}\left(x\right)u_{\theta^g}\left(x\right),\nonumber
\end{eqnarray}
\begin{eqnarray}
J_{G}&=&\displaystyle A^2 \beta^2\sum_{x=0}^{\infty} e^{\beta f^A_{\theta^g}\left(x\right)}f^{A}_{\theta^g}\left(x\right)\left(2f^A_{\theta^g}\left(x\right)-g^A\left(x\right)\right)u_{\theta^g}\left(x\right)u^{T}_{\theta^g}\left(x\right)\nonumber\\&+&A^3 \beta^3 \sum_{x=0}^{\infty}e^{\beta f^A_{\theta^g}\left(x\right)}f^{2A}_{\theta^g}\left(x\right)\left(f^A_{\theta^g}\left(x\right)-g^A\left(x\right)\right)u_{\theta^g}\left(x\right)u^{T}_{\theta^g}\left(x\right)\nonumber\\&+&(A+B)\sum_{x=0}^{\infty}f^{B}_{\theta^g}\left(x\right)\left(\left(A+B\right)f^{A}_{\theta^g}\left(x\right)-Bg^A\left(x\right)\right)u_{\theta^g}\left(x\right)u^{T}_{\theta^g}\left(x\right)\nonumber\\&+&A^2 \beta^2\sum_{x=0}^{\infty} e^{\beta f^A_{\theta^g}\left(x\right)}f^{A}_{\theta^g}\left(x\right)\left(g^A\left(x\right)-f^A_{\theta^g}\left(x\right)\right)i_{\theta^g}\left(x\right)\nonumber\\&-&(A+B)\sum_{x=0}^{\infty}f^{B}_{\theta^g}\left(f^A_{\theta^g}\left(x\right)-g^A\left(x\right)\right)i_{\theta^g}\left(x\right).\nonumber
\end{eqnarray}
If the distribution $G$ belongs to the model family $\mathcal{F}$ with $g=f_{\theta}$, then the influence function reduces to,
\begin{eqnarray}
IF\left(y,F_{\theta},T_{\alpha,\lambda,\beta}\right)&=&J_{F_{\theta}}^{-1}N_{F_{\theta}}\left(y\right), \rm{where},\\
%\end{eqnarray}
%\begin{eqnarray}
\label{ifex}
J_{F_{\theta}}&=&\displaystyle \sum_{x=0}^{\infty}\left(A^2 \beta^2 e^{\beta f^A_{\theta}\left(x\right)}f^{2A}_{\theta}\left(x\right)+(A+B)f^{A+B}_{\theta}\left(x\right)\right)u_{\theta}\left(x\right)u^{T}_{\theta}\left(x\right),\nonumber\\
N_{F_{\theta}}\left(y\right)&=&\displaystyle A^2 \beta^2 e^{\beta f^A_{\theta}\left(y\right)}f^{2A-1}_{\theta}\left(y\right)u_{\theta}\left(y\right)+(A+B)f^{A+B-1}_{\theta}\left(y\right)u_{\theta}\left(y\right)\nonumber\\&-&\sum_{x=0}^{\infty} A^2 \beta^2 e^{\beta f^A_{\theta}\left(x\right)}f^{2A}_{\theta}\left(x\right)u_{\theta}\left(x\right)-\sum_{x=0}^{\infty}(A+B)f^{A+B}_{\theta}\left(x\right)u_{\theta}\left(x\right)\nonumber.
\end{eqnarray}
Evidently, the influence function is dependent on all the three tuning parameters. Whenever the matrix $J_{F_\theta}$ is non singular, the boundedness of the influence function depends on the ability of the coefficients to control the score function $u_\theta(y)$ in the first two terms of the numerator. In most parametric models including all exponential family models, $f_\theta^\tau(y) u_\theta(y)$ remains bounded for any $\tau > 0$; in the case $\tau = 0$, however the expression equals $u_\theta(y)$ and there is no control over it to keep it bounded. For the second term of the numerator in Equation (\ref{ifex}), this is achieved when $A + B > 1$, i.e., when $\alpha > 0$. The first term of the numerator contains an additional exponential term. However, given that $f_\theta(y) \leq 1$ for any value $y$ in the support of a discrete random variable, the first term of the numerator is easily seen to be bounded for any fixed non-zero real $\beta$ when $2A - 1 > 0$, i.e., $A > 1/2$. We now list the different possible cases for boundedness of the influence function as follows: 

\begin{enumerate}
	
	\item $\beta = 0$; here the first and third terms of the numerator vanish, and the only other condition necessary is $A + B > 1$, 
	i.e., $\alpha > 0$. This is essentially the $S$-divergence case, and shows that all minimum $S$-divergence functionals with $
	\alpha > 0$ have bounded influence (irrespective of the value of $\lambda$). In this case the allowable region for the triplet $(\alpha, \lambda, \beta)$ for bounded influence is ${\mathbb S}_1 = \left(\alpha > 0, \lambda \in {\mathbb R}, \beta = 0\right)$.  
	
	\item $\beta \neq 0$, $A = 0$. In this case also  the first and third terms of the numerator drop out and the additional required condition is $\alpha > 0$. However, since $A = 1 + \lambda (1 - \alpha) = 0$, this implies $\lambda = - \frac{1}{1-\alpha}$. In this case the influence function is independent of $\beta$.  Now the relevant region for the triplet is 
	${\mathbb S}_2 = \left(\alpha>0, \lambda = - \frac{1}{1-\alpha}, \beta \neq 0\right)$.
	
	\item Now suppose $A + B = 0$, without the components being individually zero. In this case the second and fourth terms get eliminated 
	and we have $\alpha = -1$.  In this case the condition $2A - 1 > 0$ translates to $\lambda > -\frac{1}{4}$. Here the corresponding region for the triplet is ${\mathbb S}_3 = \left(\alpha = -1, \lambda \geq  - \frac{1}{4}, \beta \neq 0\right)$.
	
	\item Now we allow all the terms $\beta$, $A$ and $A+B$ to be non-zero. In this case all the four terms of the numerator are non-vanishing.
	Then, beyond the condition on $\beta$, the required conditions are $\alpha > 0$ and $\lambda\left(1 - \alpha\right) > -\frac{1}{2}$. The region here is 
	${\mathbb S}_4 = \left(\alpha > 0, \lambda (1 - \alpha) > -\frac{1}{2}, \beta \neq 0\right)$.
\end{enumerate}

Combining all the cases, we see that the IF will be bounded if the triplet $\left(\alpha,\lambda,\beta\right)\in \mathbb{S}={\mathbb S}_1 \cup {\mathbb S}_2 \cup {\mathbb S}_3 \cup {\mathbb S}_4$. 

\begin{figure}[t]%\label{pic}
	\centering
	\begin{tabular}{@{}cc@{}} 
		\includegraphics[height=.22\textheight,width=0.45\textwidth]{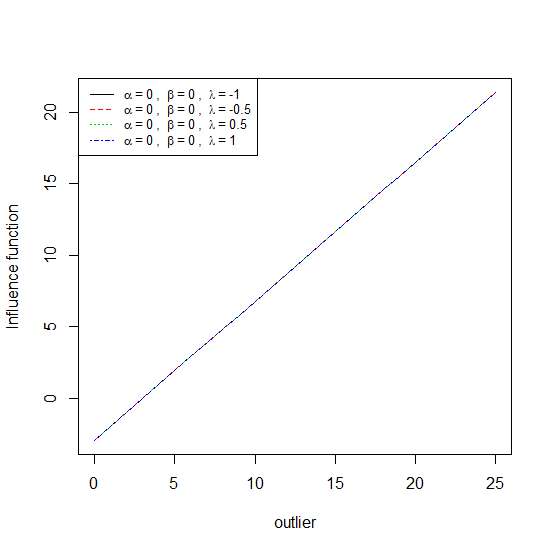}&
		\includegraphics[height=.22\textheight,width=0.45\textwidth]{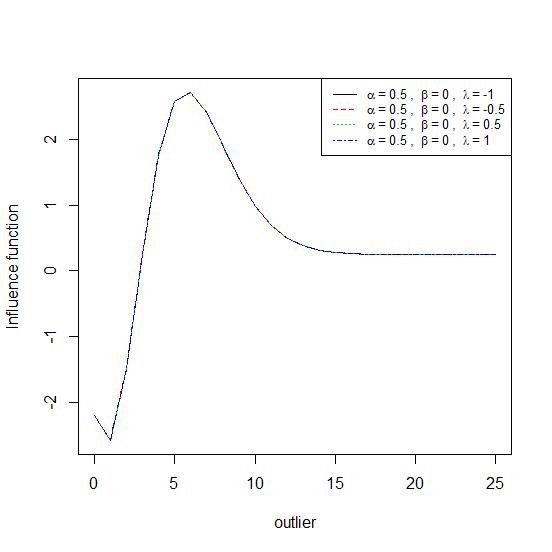}\\ 
		\includegraphics[height=.22\textheight,width=0.45\textwidth]{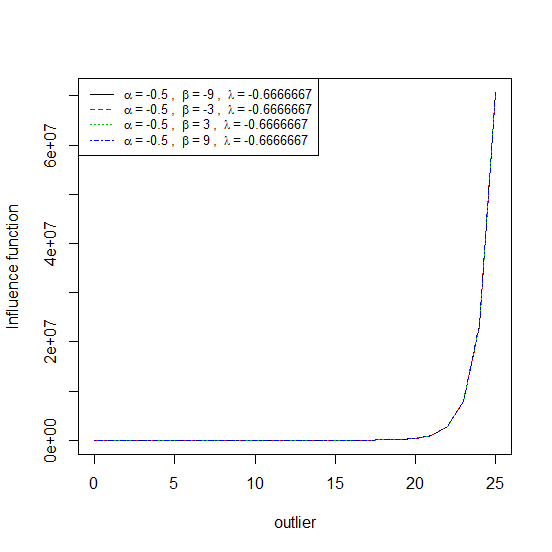}& 
		\includegraphics[height=.22\textheight,width=0.45\textwidth]{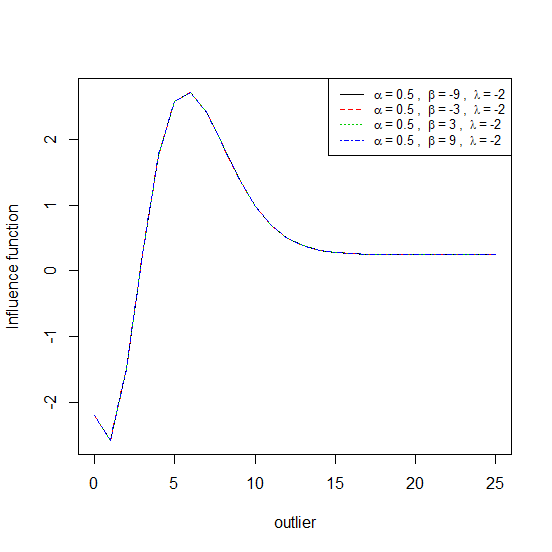}\\ 
		\includegraphics[height=.22\textheight,width=0.45\textwidth]{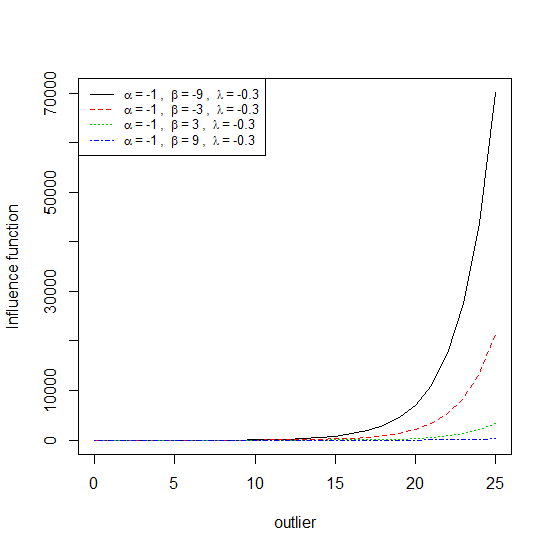} &
		\includegraphics[height=.22\textheight,width=0.45\textwidth]{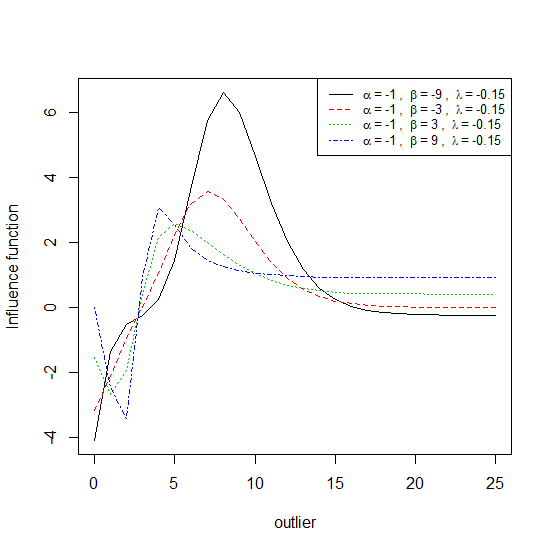}\\ 
		\includegraphics[height=.22\textheight,width=0.45\textwidth]{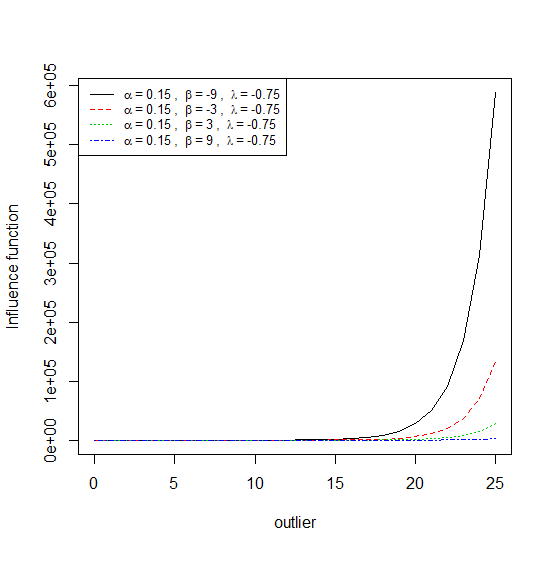}& 
		\includegraphics[height=.22\textheight,width=0.45\textwidth]{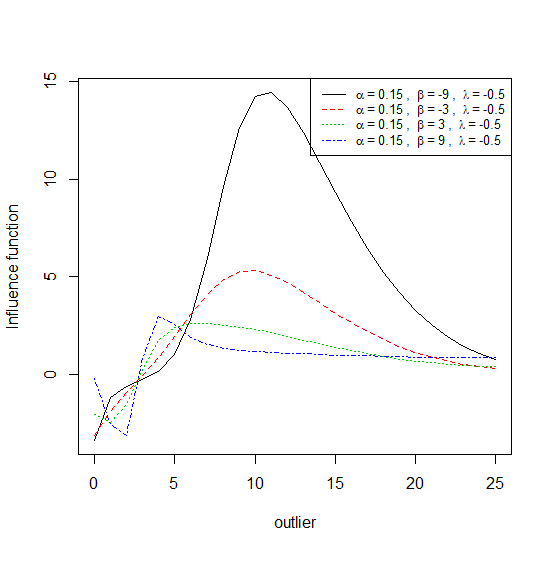}\\ 
	\end{tabular}
	\caption{Examples of unbounded influence functions (left panel) and bounded influence functions (right panel) corresponding to $\left(\alpha,\lambda,\beta\right) \in$ \rm{each disjoint subsets contained in} $\mathbb{S}$.}
	\label{f1}
\end{figure}

It is easily seen that the four constituent subregions are disjoint. For illustration, we present some plots for bounded and unbounded influence functions for the minimum GSB functional under the Poisson($\theta$) model in Figure \ref{f1}, where the true data distribution is Poisson(3). In the four rows of the right panel we give examples of triplets belonging to the four disjoint components of ${\mathbb S}$. In the first two rows of the right panel, the $\alpha$ value alone determines the shape of the curve. On the $i$-th row of the left panel, on the other hand, the triplets are slightly different from the triplets of $i$-th row on the right, but far enough to be pushed out of ${\mathbb S}_i$. Accordingly, all the plots on the left correspond to unbounded influence functions. Generally, it may also be observed that for increasing $\beta$ the curves get flatter in each plot, where IF varies over different $\beta$. We will provide further illustration of the bounded influence region of the triplet through three-dimensional graphs at the end of the simulation section.

%\clearpage

\section{Simulation results}
In the simulation section our aim is to demonstrate that by choosing non-zero values of the parameter $\beta$, we may be able to generate procedures which, in a suitable sense, improve upon the estimators that are provided by the existing standard, the class of $S$-divergences. We consider the Poisson ($\theta$) model, and choose samples of size 50 from the $(1-\epsilon) {\rm Poisson} (3) + \epsilon {\rm Poisson} (10)$ mixture, where the second component is the contaminant and $\epsilon \in [0, 1)$ is the contaminating proportion. The values 0, 0.05, 0.1 and 0.2 are considered for $\epsilon$, and at each contamination level, the samples are replicated 1000 times. The Poisson parameter is estimated in each of the 1000 replications, for each contamination level, and at each of several $(\alpha, \lambda, \beta)$ triplets considered in our study. Subsequently we construct the empirical mean square error~(MSE) against the target value of 3, for each tuning parameter triplet and each contamination level over the 1000 replications.

In case of the minimum $S$-divergence estimator, \cite{m} have empirically identified a subset of $\left(\alpha,\lambda\right)$ collections which represent good choices. According to them, the zone of ``best" estimators correspond to an elliptical subset of the tuning parameter space, with $\alpha \in [0.1, 0.6]$ and $\lambda \in [-1, -0.3]$. We hope to show that for most of the $(\alpha, \lambda)$ combinations (including the best ones) there is a corresponding better or competitive  $\left(\alpha,\lambda,\beta\right)$ combination with a non-zero $\beta$, thus providing an option which appears to perform better, at least to the extent of the findings in these simulations.

We begin with an exploration of the $S$-divergence, since this is the basis for comparison. The MSEs are presented in Table \ref{tab3} over a cross-classified grid with $\alpha$ values in \{0.1, 0.25, 0.4, 0.5, 0.6, 0.8, 1\} and $\lambda$ values in 
\{-1, -0.7, -0.5, -0.3, 0, 0.2, 0.5, 0.8, 1\}, a total of 63 cells. In each cell the empirical MSEs for $\epsilon = 0, 0.05, 0.1$ and $0.2$ are presented in a column of four elements, in that order, followed by the corresponding combination of tuning parameters $(\alpha, \lambda, \beta=0)$. We have carried the $\beta = 0$ parameter in each triplet of parameters, to indicate that the $S$-divergence is indeed a special case of the GSB divergence. It may be noted that between all the cells, there is no unique $(\alpha, \lambda)$ combination which produces an overall best result~(in terms of smallest MSE) over all the four columns (levels of contamination). 

We now expand the exploration by considering, in addition, a grid of possible non-zero $\beta$ values at each $(\alpha, \lambda)$ combination to see if the results can be improved. To be conservative about our definition of improvement, we declare the existence of a ``better'' triplet in the GSB sense if {\underline{all}} the four mean square errors corresponding to a ${(\alpha, \lambda)}$ combination within the $S$-divergence family in Table 3 are improved (reduced) by a suitable member of the GSB divergence class which is strictly outside the $S$-divergence family (corresponding to a non-zero $\beta$).
\clearpage

\clearpage

Our exploration indicates that in a large majority of the 63 cells there is a member of the GSB divergence with a non-zero $\beta$ parameter which improves (over all the four cells) the performance of the corresponding $S$-divergence estimator with the same $(\alpha, \lambda)$ combination. Interestingly it turns out that in practically all the cases where an improvement is observed it happens for a negative value of $\beta$~(it is observed to be zero in rare cases, but is never positive). A more detailed inspection indicates that in many of these cases, the improvement occurs at the value $\beta = -4$.  

In order to summarize the findings of this rather large exploration~(presented in Table \ref{tab4}) in a meaningful manner, we first note the following different cases, 
\begin{enumerate}
	\item{(First Case)}
	These are the cells where all the four mean square errors for the $S$-divergence case are reduced by the minimum GSB divergence estimator with the same values of $(\alpha, \lambda)$ and $\beta = -4$. These cells are highlighted with the blue colour in Table \ref{tab4}. (There are 18 such cells). 
	
	\item{(Second Case)}
	These are the cells where all the four MSEs for the $S$-divergence case are reduced by a  minimum GSB divergence estimator with $\beta =- 4$ but with a different $(\alpha, \lambda)$ combination than that for the corresponding cell. These cells are highlighted in red in Table \ref{tab4}. (There are 39 such cells). 
	
	\item{(Third Case)}
	These are the cells where all the four MSEs are reduced by a minimum GSB divergence estimator outside the $S$-divergence family, but with $\beta \neq -4$, and not necessarily the same $(\alpha, \lambda)$. These cells are highlighted in orange in Table \ref{tab4}. (There is one such cell). 
	
	\item{(Fourth Case)}
	These are the cells where some triplet within the minimum GSB divergence class can improve upon the three MSEs under contamination ($\epsilon = 0.05, 0.1, 0.2$) but not all the four MSEs simultaneously. While these are not ``better'' triplets in the sense described earlier in the section, the pure data MSEs (not reported here) for these triplets are close to those of the S-Divergence MSEs for these cells; in this sense these triplets are at least competitive. These cells are highlighted in green in Table \ref{tab4}. (There are three such cases). 
	
	\item{(Fifth Case)}
	These are the cells where no $(\alpha, \lambda, \beta)$ provides an improvement over the S-divergence results in the sense of any of the previous four cases~(although there are competitive alternatives). These cells remain in black in Table \ref{tab4}. There are 2 such cells.  
\end{enumerate}

On the whole, therefore, it turns out that we observe improvements in 57 out of the 63 cells in all four rows of the column of MSEs in that cell by choosing $\beta = -4$ together with the $S$-divergence parameters. Even in the handful of cases (cells) where we do not have an improvement in all the rows of the column, there generally are competitive (although not strictly better)~options within the minimum GSB divergence class with a negative value of $\beta$. 

\begin{figure}[t]
	\begin{center}
		\begin{tabular}{@{}c} 
			\includegraphics[height=.5\textheight,width=0.8\textwidth]{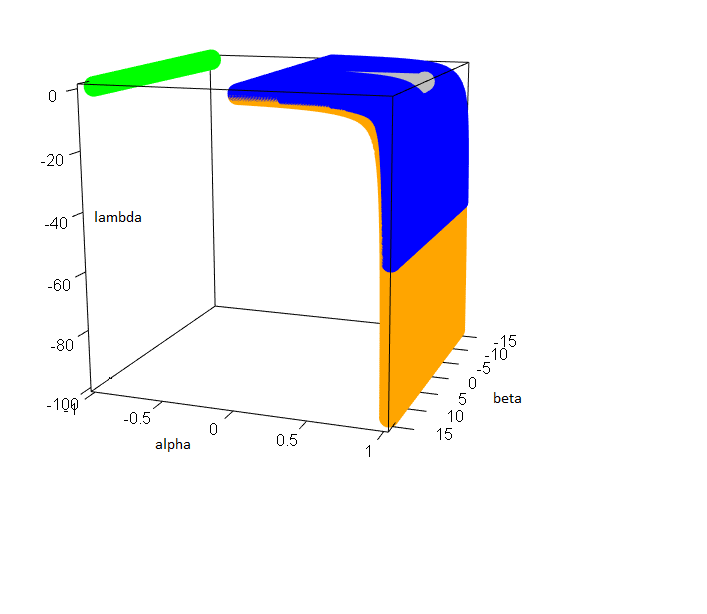}
		\end{tabular}
		\begin{tabular}{@{}c} 
			\includegraphics[height=.45\textheight,width=0.75\textwidth]{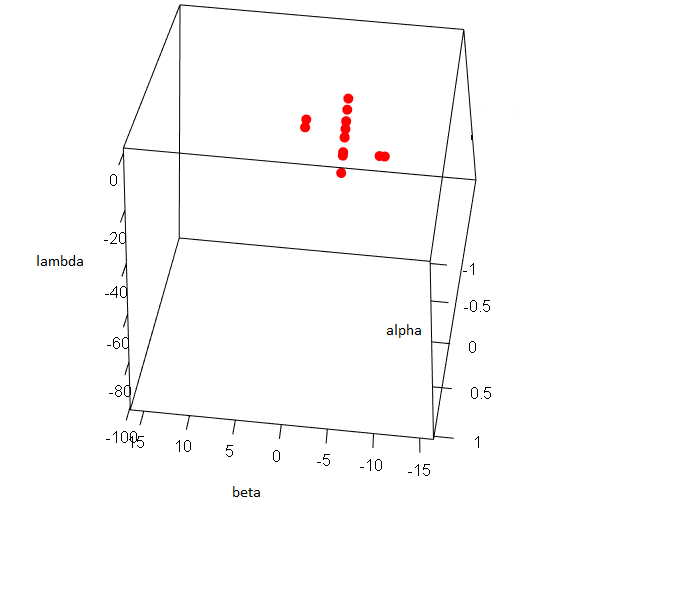}
		\end{tabular}
		\caption{The figure shows the region needed for bounded IF. Here, the grey, the orange, the green and the blue planes represent the boundaries of the sets 	${\mathbb S}_1$, ${\mathbb S}_2$, ${\mathbb S}_3$ and 	${\mathbb S}_4$, respectively.}
		\label{f2}
	\end{center}
\end{figure}

In Table \ref{tab4}, in each cell, we also present the particular $(\alpha, \lambda, \beta)$ combination which generates the mean square errors~(improved over Table \ref{tab3} in most cases, as we have seen) reported in that cell. In Figure \ref{f2}, we provide a three-dimensional plot~(as described in that section) in the three-dimensional $(\alpha, \lambda, \beta)$ plane, where the region ${\cal S}$ has been expressed as a union of several colour-coded subregions representing the individual components. The triplets corresponding to the improved MSE solutions reported in the cells of Table \ref{tab4} all belong to the blue subregion of this figure, indicating that all improved solutions are provided by bounded influence estimators.

\begin{table}[h]
	\begin{center}
		\caption{MSEs of the minimum divergence estimators within the $S$-divergence family for pure and contaminated data}
		\label{tab3}
		\begin{tabular}{c c c c c c c}
			\hline
			\hline
			0.1968 & 0.0836 & 0.0704 & 0.0708 & 0.0733 & 0.0802 & 0.0876\\
			0.1974 & 0.0981 & 0.0855 & 0.0852 & 0.0869 & 0.0926 & 0.0994\\
			0.1753 & 0.1063 & 0.1012 & 0.1028 & 0.1054 & 0.1116 & 0.1118\\
			0.3099 & 0.2245 & 0.2119 & 0.2113 & 0.2130 & 0.2200 & 0.2298\\
			(0.1, $-1$, 0) & (0.25, -1, 0) & (0.4, -1, 0) & (0.5, -1, 0) & (0.6, -1, 0) & (0.8, -1, 0) & (1, -1, 0)\\
			\hline
			0.0751 & 0.0666 & 0.0673 & 0.0698 & 0.0729 & 0.0800 & 0.0876\\
			0.0893 & 0.0830 & 0.0831 & 0.0847 & 0.0869 & 0.0927 & 0.0994\\
			0.1081 & 0.1044 & 0.1045 & 0.1056 & 0.1073 & 0.1121 & 0.1118\\
			0.2830 & 0.2505 & 0.2328 & 0.2264 & 0.2231 & 0.2233 & 0.2298\\
			(0.1, -0.7, 0) & (0.25, -0.7, 0) & (0.4, -0.7, 0) & (0.5, -0.7, 0) & (0.6, -0.7, 0) & (0.8, -0.7, 0) & (1, -0.7, 0)\\
			\hline
			0.0638 & 0.0635 & 0.0665 & 0.0694 & 0.0727 & 0.0799 & 0.0876\\
			0.0836 & 0.0821 & 0.0832 & 0.0849 & 0.0871 & 0.0927 & 0.0994\\
			0.1203 & 0.1120 & 0.1087 & 0.1083 & 0.1089 & 0.1125 & 0.1118\\
			0.3715 & 0.2958 & 0.2559 & 0.2408 & 0.2319 & 0.2258 & 0.2298\\
			(0.1, -0.5, 0) & (0.25, -0.5, 0) & (0.4, -0.5, 0) & (0.5, -0.5, 0) & (0.6, -0.5, 0) & (0.8, -0.5, 0) & (1, -0.5, 0)\\
			\hline
			0.0600 & 0.0622 & 0.0660 & 0.0691 & 0.0725 & 0.0798 & 0.0876\\
			0.0895 & 0.0846 & 0.0843 & 0.0856 & 0.0875 & 0.0928 & 0.0994\\
			0.1554 & 0.1264 & 0.1149 & 0.1118 & 0.1109 & 0.1129 & 0.1118\\
			0.5669 & 0.3709 & 0.2904 & 0.2605 & 0.2424 & 0.2286 & 0.2298\\
			(0.1, -0.3, 0) & (0.25, -0.3, 0) & (0.4, -0.3, 0) & (0.5, -0.3, 0) & (0.6, -0.3, 0) & (0.8, -0.3, 0) & (1, -0.3, 0)\\
			\hline
			0.0592 & 0.0617 & 0.0657 & 0.0688 & 0.0721 & 0.0796 & 0.0876\\
			0.1415 & 0.0971 & 0.0880 & 0.0873 & 0.0883 & 0.0929 & 0.0994\\
			0.3491 & 0.1774 & 0.1308 & 0.1196 & 0.1147 & 0.1136 & 0.1118\\
			1.3860 & 0.6555 & 0.3832 & 0.3061 & 0.2655 & 0.2333 & 0.2298\\
			(0.1, 0, 0) & (0.25, 0, 0) & (0.4, 0, 0) & (0.5, 0, 0) & (0.6, 0, 0) & (0.8, 0, 0) & (1, 0, 0)\\
			\hline
			0.0608 & 0.0621 & 0.0657 & 0.0687 & 0.0721 & 0.0794 & 0.0876\\
			0.3302 & 0.1231 & 0.0930 & 0.0892 & 0.0890 & 0.0930 & 0.0994\\
			0.8565 & 0.2745 & 0.1508 & 0.1276 & 0.1181 & 0.1141 & 0.1118\\
			2.5938 & 1.0853 & 0.5550 & 0.3553 & 0.2867 & 0.2370 & 0.2298\\
			(0.1, 0.2, 0) & (0.25, 0.2, 0) & (0.4, 0.2, 0) & (0.5, 0.2, 0) & (0.6, 0.2, 0) & (0.8, 0.2, 0) & (1, 0.2, 0)\\
			\hline
			0.0671 & 0.0638 & 0.0658 & 0.0685 & 0.0718 & 0.0792 & 0.0876\\
			1.1434 & 0.3251 & 0.1115 & 0.0943 & 0.0907 & 0.0931 & 0.0994\\
			2.3829 & 0.8165 & 0.2234 & 0.1489 & 0.1255 & 0.1151 & 0.1118\\
			4.8817 & 2.4261 & 0.8641 & 0.4847 & 0.3338 & 0.2434 & 0.2298\\
			(0.1, 0.5, 0) & (0.25, 0.5, 0) & (0.4, 0.5, 0) & (0.5, 0.5, 0) & (0.6, 0.5, 0) & (0.8, 0.5, 0) & (1, 0.5, 0)\\
			\hline
			0.0778 & 0.0676 & 0.0665 & 0.0685 & 0.0716 & 0.0790 & 0.0876\\
			1.9928 & 0.9339 & 0.1951 & 0.1068 & 0.0936 & 0.0933 & 0.0994\\
			3.7890 & 1.9909 & 0.4869 & 0.1994 & 0.1378 & 0.1162 & 0.1118\\
			6.6731 & 4.2592 & 1.6784 & 0.7520 & 0.4130 & 0.2511 & 0.2298\\
			(0.1, 0.8, 0) & (0.25, 0.8, 0) & (0.4, 0.8, 0) & (0.5, 0.8, 0) & (0.6, 0.8, 0) & (0.8, 0.8, 0) & (1, 0.8, 0)\\
			\hline
			0.0863 & 0.0717 & 0.0673 & $0.0686$ & 0.0714 & 0.0789 & 0.0876\\
			2.4449 & 1.3987 & 0.3803 & $0.1283$ & 0.0967 & 0.0934 & 0.0994\\
			4.5000 & 2.7992 & 0.9117 & $0.2793$ & 0.1514 & 0.1171 & 0.1118\\
			7.5320 & 5.3554 & 2.5215 & $1.0745$ & 0.4969 & 0.2572 & 0.2298\\
			(0.1, 1, 0) & (0.25, 1, 0) & (0.4, 1, 0) & (0.5, 1, 0) & (0.6, 1, 0) & (0.8, 1, 0) & (1, 1, 0)\\
			\hline
			\hline
		\end{tabular}
	\end{center}
\end{table}

\begin{table}[h]
	\begin{center}
		\caption{MSEs of the minimum GSB divergence estimators under pure and contaminated data}
		\label{tab4}
		\begin{tabular}{c c c c c c c}
			\hline
			\hline
			\textcolor{red}{0.0623} & \textcolor{red}{0.0696} & 0.0704 & 0.0708 & \textcolor{YellowOrange}{0.0687} & \textcolor{red}{0.0720} & \textcolor{blue}{0.0763}\\
			\textcolor{red}{0.0816} & \textcolor{red}{0.0843} & 0.0855 & 0.0852 & \textcolor{YellowOrange}{0.0833} & \textcolor{red}{0.0859} & \textcolor{blue}{0.0892}\\
			\textcolor{red}{0.1115} & \textcolor{red}{0.1056} & 0.1012 & 0.1028 & \textcolor{YellowOrange}{0.1042} & \textcolor{red}{0.1060} & \textcolor{blue}{0.1077}\\
			\textcolor{red}{0.2831} & \textcolor{red}{0.2207} & 0.2119 & 0.2113 & \textcolor{YellowOrange}{0.2110} & \textcolor{red}{0.2162} & \textcolor{blue}{0.2115}\\
			\textcolor{red}{(0.4, $-0.4$, -4)} & \textcolor{red}{(0.8, -0.5, -4)} & (0.4, -1, 0) & (0.5, -1, 0) & \textcolor{YellowOrange}{(0.8, 0, -7.5)} & \textcolor{red}{(0.8, -0.3, -4)}  & \textcolor{blue}{(1, -1, -4)}\\
			\hline
			\textcolor{red}{0.0642} & \textcolor{teal}{0.0681} & \textcolor{teal}{0.0681} & \textcolor{red}{0.0696} & \textcolor{red}{0.0696} & \textcolor{red}{0.0720} & \textcolor{blue}{0.0763}\\
			\textcolor{red}{0.0816} & \textcolor{teal}{0.0826} & \textcolor{teal}{0.0826} & \textcolor{red}{0.0843} & \textcolor{red}{0.0843} & \textcolor{red}{0.0859} & \textcolor{blue}{0.0892}\\
			\textcolor{red}{0.1076} & \textcolor{teal}{0.1043} & \textcolor{teal}{0.1043} & \textcolor{red}{0.1055} & \textcolor{red}{0.1055} & \textcolor{red}{0.1060} & \textcolor{blue}{0.1077}\\
			\textcolor{red}{0.2514} & \textcolor{teal}{0.2135} & \textcolor{teal}{0.2135} & \textcolor{red}{0.2207} & \textcolor{red}{0.2207} & \textcolor{red}{0.2162} & \textcolor{blue}{0.2115}\\
			\textcolor{red}{(0.6, -0.5, -4)} & \textcolor{teal}{(0.8, 0, -8)} & \textcolor{teal}{(0.8, 0, -8)} & \textcolor{red}{(0.8, -0.5, -4)} & \textcolor{red}{(0.8, -0.5, -4)} & \textcolor{red}{(0.8, -0.3, -4)}  & \textcolor{blue}{(1, -0.7, -4)}\\
			\hline
			\textcolor{red}{0.0623} & \textcolor{red}{0.0623} & \textcolor{red}{0.0659} & \textcolor{red}{0.0678} & \textcolor{red}{0.0678} & \textcolor{blue}{0.0696} & \textcolor{blue}{0.0763}\\
			\textcolor{red}{0.0816} & \textcolor{red}{0.0816} & \textcolor{red}{0.0825} & \textcolor{red}{0.0834} & \textcolor{red}{0.0834} & \textcolor{blue}{0.0843} & \textcolor{blue}{0.0892}\\
			\textcolor{red}{0.1115} & \textcolor{red}{0.1115} & \textcolor{red}{0.1071} & \textcolor{red}{0.1061} & \textcolor{red}{0.1061} & \textcolor{blue}{0.1055} & \textcolor{blue}{0.1077}\\
			\textcolor{red}{0.2831} & \textcolor{red}{0.2831} & \textcolor{red}{0.2417} & \textcolor{red}{0.2295} & \textcolor{red}{0.2295} & \textcolor{blue}{0.2207} & \textcolor{blue}{0.2115}\\
			\textcolor{red}{(0.4, -0.4, -4)} & \textcolor{red}{(0.4, -0.4, -4)} & \textcolor{red}{(0.5, -0.3, -4)} & \textcolor{red}{(0.6, -0.3, -4)} & \textcolor{red}{(0.6, -0.3, -4)} & \textcolor{blue}{(0.8, -0.5, -4)}  & \textcolor{blue}{(1, -0.5, -4)}\\
			\hline
			\textcolor{blue}{0.0600} & \textcolor{blue}{0.0619} & \textcolor{blue}{0.0642} & \textcolor{blue}{0.0659} & \textcolor{blue}{0.0678} & \textcolor{blue}{0.0720} & \textcolor{blue}{0.0763}\\
			\textcolor{blue}{0.0845} & \textcolor{blue}{0.0822} & \textcolor{blue}{0.0819} & \textcolor{blue}{0.0825} & \textcolor{blue}{0.0834} & \textcolor{blue}{0.0859} & \textcolor{blue}{0.0892}\\
			\textcolor{blue}{0.1294} & \textcolor{blue}{0.1154} & \textcolor{blue}{0.1091} & \textcolor{blue}{0.1071} & \textcolor{blue}{0.1061} & \textcolor{blue}{0.1060} & \textcolor{blue}{0.1077}\\
			\textcolor{blue}{0.4049} & \textcolor{blue}{0.3080} & \textcolor{blue}{0.2602} & \textcolor{blue}{0.2417} & \textcolor{blue}{0.2295} & \textcolor{blue}{0.2162} & \textcolor{blue}{0.2115}\\
			\textcolor{blue}{(0.1, -0.3, -4)} & \textcolor{blue}{(0.25, -0.3, -4)} & \textcolor{blue}{(0.4, -0.3, -4)} & \textcolor{blue}{(0.5, -0.3, -4)} & \textcolor{blue}{(0.6, -0.3, -4)} & \textcolor{blue}{(0.8, -0.3, -4)}  & \textcolor{blue}{(1, -0.3, -4)}\\
			\hline
			\textcolor{teal}{0.0600} & \textcolor{red}{0.0600} & \textcolor{red}{0.0644} & \textcolor{red}{0.0644} & \textcolor{red}{0.0644} & \textcolor{blue}{0.0755} & \textcolor{blue}{0.0763}\\
			\textcolor{teal}{0.0845} & \textcolor{red}{0.0845} & \textcolor{red}{0.0817} & \textcolor{red}{0.0817} & \textcolor{red}{0.0817} & \textcolor{blue}{0.0887} & \textcolor{blue}{0.0892}\\
			\textcolor{teal}{0.1294} & \textcolor{red}{0.1294} & \textcolor{red}{0.1073} & \textcolor{red}{0.1073} & \textcolor{red}{0.1073} & \textcolor{blue}{0.1077} & \textcolor{blue}{0.1077}\\
			\textcolor{teal}{0.4049} & \textcolor{red}{0.4049} & \textcolor{red}{0.2492} & \textcolor{red}{0.2492} & \textcolor{red}{0.2492} & \textcolor{blue}{0.2139} & \textcolor{blue}{0.2115}\\
			\textcolor{teal}{(0.1, -0.3, -4)} & \textcolor{red}{(0.1, -0.3, -4)} & \textcolor{red}{(0.8, -1, -4)} & \textcolor{red}{(0.8, -1, -4)} & \textcolor{red}{(0.8, -1, -4)} & \textcolor{blue}{(0.8, 0, -4)} & \textcolor{blue}{(1, 0, -4)}\\
			\hline
			\textcolor{red}{0.0600} & \textcolor{red}{0.0600} & \textcolor{red}{0.0644} & \textcolor{red}{0.0644} & \textcolor{red}{0.0644} & \textcolor{blue}{0.0779} & \textcolor{blue}{0.0763}\\
			\textcolor{red}{0.0845} & \textcolor{red}{0.0845} & \textcolor{red}{0.0817} & \textcolor{red}{0.0817} & \textcolor{red}{0.0817} & \textcolor{blue}{0.0907} & \textcolor{blue}{0.0892}\\
			\textcolor{red}{0.1294} & \textcolor{red}{0.1294} & \textcolor{red}{0.1073} & \textcolor{red}{0.1073} & \textcolor{red}{0.1073} & \textcolor{blue}{0.1092} & \textcolor{blue}{0.1077}\\
			\textcolor{red}{0.4049} & \textcolor{red}{0.4049} & \textcolor{red}{0.2492} & \textcolor{red}{0.2492} & \textcolor{red}{0.2492} & \textcolor{blue}{0.2146} & \textcolor{blue}{0.2115}\\
			\textcolor{red}{(0.1, -0.3, -4)} & \textcolor{red}{(0.1, -0.3, -4)} & \textcolor{red}{(0.8, -1, -4)} & \textcolor{red}{(0.8, -1, -4)} & \textcolor{red}{(0.8, -1, -4)} & \textcolor{blue}{(0.8, 0.2, -4)} & \textcolor{blue}{(1, 0.2, -4)}\\
			\hline
			\textcolor{red}{0.0600} & \textcolor{red}{0.0600} & \textcolor{red}{0.0644} & \textcolor{red}{0.0644} & \textcolor{red}{0.0644} & \textcolor{red}{0.0696} & \textcolor{blue}{0.0763}\\
			\textcolor{red}{0.0845} & \textcolor{red}{0.0845} & \textcolor{red}{0.0817} & \textcolor{red}{0.0817} & \textcolor{red}{0.0817} & \textcolor{red}{0.0843} & \textcolor{blue}{0.0892}\\
			\textcolor{red}{0.1294} & \textcolor{red}{0.1294} & \textcolor{red}{0.1073} & \textcolor{red}{0.1073} & \textcolor{red}{0.1073} & \textcolor{red}{0.1055} & \textcolor{blue}{0.1077}\\
			\textcolor{red}{0.4049} & \textcolor{red}{0.4049} & \textcolor{red}{0.2492} & \textcolor{red}{0.2492} & \textcolor{red}{0.2492} & \textcolor{red}{0.2207} & \textcolor{blue}{0.2115}\\
			\textcolor{red}{(0.1, -0.3, -4)} & \textcolor{red}{(0.1, -0.3, -4)} & \textcolor{red}{(0.8, -1, -4)} & \textcolor{red}{(0.8, -1, -4)} & \textcolor{red}{(0.8, -1, -4)} & \textcolor{red}{(0.8, -0.5, -4)} & \textcolor{blue}{(1, 0.5, -4)}\\
			\hline
			\textcolor{red}{0.0600} & \textcolor{red}{0.0600} & \textcolor{red}{0.0644} & \textcolor{red}{0.0644} & \textcolor{red}{0.0644} & \textcolor{red}{0.0696} & \textcolor{blue}{0.0763}\\
			\textcolor{red}{0.0845} & \textcolor{red}{0.0845} & \textcolor{red}{0.0817} & \textcolor{red}{0.0817} & \textcolor{red}{0.0817} & \textcolor{red}{0.0843} & \textcolor{blue}{0.0892}\\
			\textcolor{red}{0.1294} & \textcolor{red}{0.1294} & \textcolor{red}{0.1073} & \textcolor{red}{0.1073} & \textcolor{red}{0.1073} & \textcolor{red}{0.1055} & \textcolor{blue}{0.1077}\\
			\textcolor{red}{0.4049} & \textcolor{red}{0.4049} & \textcolor{red}{0.2492} & \textcolor{red}{0.2492} & \textcolor{red}{0.2492} & \textcolor{red}{0.2207} & \textcolor{blue}{0.2115}\\
			\textcolor{red}{(0.1, -0.3, -4)} & \textcolor{red}{(0.1, -0.3, -4)} & \textcolor{red}{(0.8, -1, -4)} & \textcolor{red}{(0.8, -1, -4)} & \textcolor{red}{(0.8, -1, -4)} & \textcolor{red}{(0.8, -0.5, -4)} & \textcolor{blue}{(1, 0.8, -4)}\\
			\hline
			\textcolor{red}{0.0600} & \textcolor{red}{0.0600} & \textcolor{red}{0.0644} & \textcolor{red}{0.0644} & \textcolor{red}{0.0644} & \textcolor{red}{0.0696} & \textcolor{blue}{0.0763}\\
			\textcolor{red}{0.0845} & \textcolor{red}{0.0845} & \textcolor{red}{0.0817} & \textcolor{red}{0.0817} & \textcolor{red}{0.0817} & \textcolor{red}{0.0843} & \textcolor{blue}{0.0892}\\
			\textcolor{red}{0.1294} & \textcolor{red}{0.1294} & \textcolor{red}{0.1073} & \textcolor{red}{0.1073} & \textcolor{red}{0.1073} & \textcolor{red}{0.1055} & \textcolor{blue}{0.1077}\\
			\textcolor{red}{0.4049} & \textcolor{red}{0.4049} & \textcolor{red}{0.2492} & \textcolor{red}{0.2492} & \textcolor{red}{0.2492} & \textcolor{red}{0.2207} & \textcolor{blue}{0.2115}\\
			\textcolor{red}{(0.1, -0.3, -4)} & \textcolor{red}{(0.1, -0.3, -4)} & \textcolor{red}{(0.8, -1, -4)} & \textcolor{red}{(0.8, -1, -4)} & \textcolor{red}{(0.8, -1, -4)} & \textcolor{red}{(0.8, -0.5, -4)} & \textcolor{blue}{(1, 1, -4)}\\
			\hline
			\hline
		\end{tabular}
	\end{center}
\end{table}
\clearpage

\clearpage

\section{Selection of tuning parameters}
Our simulations in the previous section seem to suggest that the minimum divergence estimators within the GSB class with $\beta = -4$ often provides good options for data analysis. To take full advantage of this observation, this subclass of the GSB family should be explored further. However, we want to fully exploit the flexibility of the three parameter system, and noting that in some cases the optimal is outside the $\beta = -4$ subclass, including some which generate the most competitive solutions in the full system, we want to use an overall 
data-based tuning parameter selection rule in which all the three parameters are allowed to vary over reasonable supports. The aim is to select the ``best'' tuning parameter combination depending on the amount of anomaly in the data. Thus datasets which show very close compatibility to the model should be analyzed by a triplet providing an efficient solution, while a more anomalous one should have a more robust member of the GSB class to deal with it. 

The current literature contains some suggestions for choosing data-driven tuning parameters in specific situations which we can make use of. The works of \cite{r}, \cite{s} and \cite{t} are relevant in this connection. Different algorithms for generating the ``optimal" tuning parameter in case of the DPD have been described in these papers, which we will denote as the HK (for Hong and Kim), the OWJ (for one-step Warwick-Jones) and the IWJ (for the iterated Warwick-Jones) algorithms, respectively. The essential idea here is to construct an empirical approximation to the mean square error as a function of the tuning parameters (and a pilot estimator) and minimize it over the tuning parameter. The IWJ algorithm of \cite{t} refines the OWJ algorithm of \cite{s} by choosing the solution at a particular stage as the pilot for the next stage and going through an iterative process, thus eliminating the dependence on the pilot estimator as long as one has a good robust estimator to start with. The Hong and Kim algorithm does not consider the bias part of the mean square error and occasionally throws up highly non-robust solutions. See \cite{t} for a full description  and a comparative discussion of all the three algorithms. We will implement the same algorithms here, but for the GSB parameters rather than just the DPD parameter. 

In the following we have taken up two real data examples and considered the problem of selecting the ``optimal" tuning parameters in each case. The OWJ algorithm considered here uses the minimum $L_2$ distance estimator as the pilot. Although the IWJ algorithm is pilot independent, for computational purposes it needs to commence from some suitable robust pilot for which also we utilize the minimum $L_2$ distance estimator. While the IWJ algorithm is our preferred method, we demonstrate the use of all the three algorithms in the following data sets. 

\begin{example} \label{peri} (Peritonitis Data): This example involves the incidence of 
	peritonitis in $390$ kidney patients. The data are available in Chapter 2 of \cite{f}, and are also presented in the Online Supplement. The observations at 10 and 12 may be regarded as mild outliers. A geometric model with success probability $\theta$ has been fitted to these frequency data. Here, the IWJ solution coincides with the HK solution where the estimate of success probability is 0.5110 corresponding to $\left(\alpha,\lambda,\beta\right)=\left(0.41,-0.84,-3.5\right)$. The OWJ solution gives a slightly different success probability of 0.5105 corresponding to  $\left(\alpha,\lambda,\beta\right)=\left(0.17,-0.60,-3\right)$. In case of clean data these IWJ, OWJ and HK estimates will be 0.5044, 0.5061 and 0.5029 corresponding to $\left(\alpha,\lambda,\beta\right)=\left(0.47,-1,-2\right)$, $\left(0.29,-1,-1\right)$ and $\left(0.55,-1,-3\right)$,  respectively, being slightly different from each other. On the contrary, the MLEs for the full dataset and the (two) outlier 
	deleted dataset are $0.4962$ and $0.5092$, respectively.
\end{example}

Now we consider a more recent dataset for the implementation of our new proposal.

\begin{figure}[ht!]
\centering
\begin{tabular}{@{}c} 
	\includegraphics[height=.4\textheight,width=0.7\textwidth]{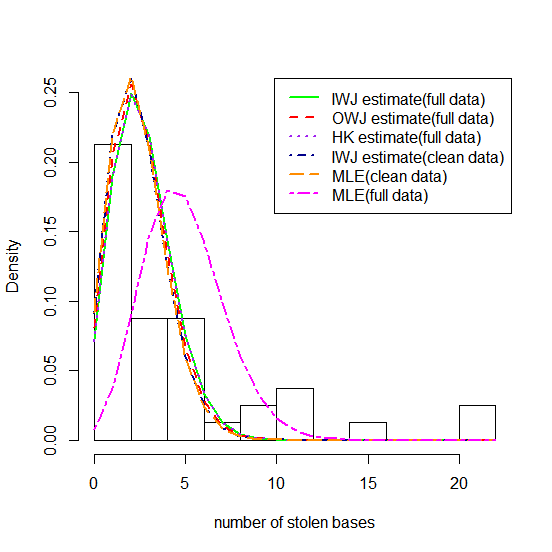}\\
\end{tabular}
\caption{Some significant fits for the Stolen Bases Data under the Poisson model. Here ``clean data" refer to the modified data after removing all 9 outliers.}
\label{f3}
\end{figure}

\begin{example} \label{stbs} (Stolen Bases Data): In ``Major League Baseball (MLB) Player Batting Stats" for the 2019 MLB Regular Season, obtained from the ESPN.com website, one  variable of interest is the number of Stolen Bases (SB) awarded to the top 40 Home Run (HR) scorers of the American League (AL). This dataset, containing three extreme and six moderate outliers, could be well-modelled by the Poisson distribution if not for the outliers. We are interested to estimate $\theta$, the average number of Stolen Bases (SB) awarded to the MLB batters of the AL throughout the whole regular season. The ``optimal'' estimates, derived from the implementation of the three algorithms under the Poisson model, are presented in Table \ref{tab5}. The fitted curves corresponding to some of these optimal estimates are given in Figure \ref{f3}. It is clear that except for the full data MLE, all the other estimators primarily model the main model conforming part of the data and sacrifice the outliers. 
\end{example}

	\begin{table}[h]
	\begin{center}
		\caption{Optimal estimates in different cases for the Stolen Bases Data}
		\begin{threeparttable}
			\label{tab5}
			\begin{tabular}{c| c| c| c}
				\hline
				data & method & optimal $\hat{\theta}$ & optimal $\left(\alpha, \lambda, \beta\right)$\\
				\hline
				Full data (with outliers) & IWJ & $2.6270$ & $\left(0.65,-0.98,-8\right)$\\
				& OWJ & $2.5086$ & $\left(0.73,-1,-8\right)$\\
				& HK & $2.6409$ & $\left(0.65,-1,-8\right)$\\
				& MLE & $4.875$ & $\left(0,0,0\right)$\\
				\hline
				excluding 9 outliers & IWJ & $2.3949$ & $\left(0.01,1.00,0\right)$\\& OWJ & $2.3229$ & $\left(0.25,1.00,0\right)$\\
				& HK & $2.6918$ & $\left(0.45,-1,-8\right)$\\
				& MLE & $2.3871$ & $\left(0,0,0\right)$\\
				\hline
			\end{tabular}
		\end{threeparttable}
	\end{center}
\end{table}

\clearpage

\section{Concluding remarks}
In this paper, we have provided an extension of the ordinary Bregman divergence which has direct applications to developing new classes of divergence measures, and, in turn, in providing more options for minimum distance inference with better mixes of model efficiency and robustness. In the second part of the paper, we have made use of the suggested approach in generating a particular super-family of divergences which seems to work very well in practice and provides new minimum divergence techniques which appear to improve the performance of the $S$-divergence based procedures in many cases. Since the results presented here are based on a single study, more research will certainly be needed to decide to what extent the observed advantages of the procedures considered here can be generalized, but clearly there appears to be enough evidence to suggest such explorations are warranted. 

Even apart from the search for other divergences, several possible follow ups of this research immediately present themselves. This paper is restricted to discrete parametric models. An obvious follow up step is to suitably handle the case of continuous models, where the construction of the density and the divergence are more difficult. Another obvious extension will be to extend the procedures to more complicated data structures beyond the simple independently and identically distributed data scenarios. Yet another extension would be to apply this and other similarly developed super-divergences in the area of robust testing of hypothesis. We hope to take up all of these extension in our future work.   

The subfamily of the class of GSB divergences with $\beta = -4$ also needs some attention, and we hope to take it up in the future. For the time being we have presented the results for our real data examples for the $\beta = -4$ subfamily of GSB in the Online Supplement.

%\section*{Acknowledgements}
%I would like to thank my PhD supervisor, Prof. Ayanendranath Basu, who is also my co-author here, for helping me immensely in preparing my paper. I would also like to thank the anonymous reviewers and the associate editor for their helpful inputs and suggestions which helped me a lot to improve the quality of the paper.

\section*{Disclosure statement}
No potential conflict of interest was reported by the authors.

\end{document}